\documentclass[a4paper,11pt]{article}
\usepackage{url}
\usepackage[utf8]{inputenc}
\usepackage{amsmath,accents,amssymb,amsthm}
\usepackage[usenames]{color}
\usepackage{comment}

\theoremstyle{plain}
\newtheorem{theorem}{Theorem}[section]
\newtheorem{lemma}[theorem]{Lemma}

\theoremstyle{definition}

\newtheorem*{example}{Example}

\newcommand\xch[2]{#1}
\newcommand\reftext[1]{#1}
\newcommand\setitemindent[1]{}
\newcommand\setleftskip[1]{}
\newcommand\tvspace[1]{#1}

\usepackage[mathscr]{eucal}
\usepackage[usenames]{color}

\usepackage{algorithm,algpseudocode}

\algnewcommand{\Or}{\textbf{ or }}
\algnewcommand{\And}{\textbf{ and }}
\algnewcommand{\From}{\textbf{ from }}
\algnewcommand{\To}{\textbf{ to }}
\algnewcommand{\Let}{\textbf{Let }}
\algnewcommand{\Var}{\texttt}
\algnewcommand{\Int}{\textbf{ int }}
\algnewcommand{\List}{\textbf{ list }}
\algnewcommand{\Bool}{\textbf{ bool }}
\algnewcommand{\True}{\textsc{ true }}
\algnewcommand{\False}{\textsc{ false }}
\algnewcommand{\nlind}{\newline\phantom{mmmmmmmmmmmmmmm}}
\renewcommand\leq\leqslant
\renewcommand\geq\geqslant

\newcommand\shrink[1]{\overline{#1}}

\usepackage{tikz,tikz-qtree,tikz-qtree-compat,adjustbox,xcolor,array,hhline}
\usepackage{pst-plot}
\usepackage{etoolbox}
\usetikzlibrary{positioning}
\providecommand\circledcolorednumb{}\renewcommand\circledcolorednumb[2]{\resizebox{0.101111\textwidth}{!}{\tikz[baseline=(char.center)]{\node[shape = circle,draw, inner sep = 2pt,fill=#1](char)    {\phantom{00}};\node[anchor=center] at (char.center) {\makebox(0,0){\large{{\sf #2}}}};}}}

\robustify{\circledcolorednumb}

\newcommand\mut[1]{\ignorespaces}

\title{Exploring the unleaved tree of numerical semigroups\\up to a given genus}
\author{Maria Bras-Amorós\\Universitat Politècnica de Catalunya}

\begin{document}
\maketitle

\begin{abstract}
We present a new algorithm for exploring or counting the numerical semigroups
of a given genus that uses the unleaved version of the tree of numerical
semigroups. In the unleaved tree there are no leaves rather than the ones
at depth equal to the genus in consideration. For exploring the unleaved
tree we present a new encoding system of a numerical semigroup given by
the gcd of its left elements and its shrinking, that is, the semigroup
generated by its left elements divided by their gcd. We show a method to
determine the right generators and strong generators of a semigroup by
means of the gcd and the shrinking encoding, as well as a method to encode
a semigroup from the encoding of its parent or of its preceding sibling.
With the new algorithm, we obtained $n_{76}=29028294421710227$ and
$n_{77}=47008818196495180$.
\end{abstract}

\begin{center}
\end{center}

%
\section{Introduction}
\label{sec1}

A \emph{numerical semigroup} is an arithmetic object consisting of a subset
$\Lambda $ of ${\mathbb N}$ that contains $0$, is closed under addition,
and has a finite complement in ${\mathbb N}_{0}$. The elements in
${\mathbb N}_{0}\setminus \Lambda $ are called the \emph{gaps} of the semigroup,
the largest gap is called the \emph{Frobenius number}, and the number of
gaps is the \emph{genus} $g(\Lambda )$ of the numerical semigroup.

There have been many efforts to compute the sequence $n_{g}$ counting the
number of numerical semigroups of genus $g$ and today the sequence values
are known up to $n_{75}$
\cite{NivaldoMedeiros,Komeda,Br:Fibonacci,FromentinHivert,rgd,DelgadoEliahouFromentin}.
See \cite{Sloane} (\texttt{\url{https://oeis.org/A007323}}) for the complete list
and for more information. In 2007 it was conjectured that the sequence
$n_{g}$ is increasing, that each term is at least the sum of the two previous
terms, and that the ratio between each term and the sum of the two previous
terms approaches one as $g$ grows to infinity, which is equivalent to having
a growth rate approaching the golden ratio $\frac{1+\sqrt{5}}{2}$
\cite{Segovia,Br:Fibonacci}. The last statement of the conjecture was proved
\xch{by}{by Alex Zhai} \cite{Zhai}.

A numerical semigroup $\Lambda $ is \emph{generated} by a set of integers
$a_{1},\dots ,a_{k}$ if
$\Lambda =a_{1}{\mathbb N}_{0}+\dots +a_{k}{\mathbb N}_{0}$. In this case
$a_{1},\dots ,a_{k}$ need to be coprime and we write
$\Lambda =\langle a_{1},\dots ,a_{k}\rangle $. Each numerical semigroup
has a unique minimal set of generators. The set of elements smaller than
the Frobenius number is called the \emph{left elements}. The minimal generators
that are not left elements are called \emph{right generators}. The elements
of $\Lambda $ can be enumerated in increasing order as
$\lambda _{0}=0,\lambda _{1},\dots $.

Special cases of numerical semigroups are \emph{ordinary} semigroups, whose
gaps are all consecutive from the integer $1$, \emph{quasi-ordinary} semigroups,
whose gaps are all consecutive from the integer $1$, except for an isolated
gap, indeed, the Frobenius number, and \emph{pseudo-ordinary} semigroups,
which are the union of an ordinary semigroup and an element not in the
ordinary semigroup.

The problem of defining an algorithm to count or visit all numerical semigroups
up to a given genus has commonly been tackled using the so-called \emph{tree
of numerical semigroups} ${\mathscr T}$. This tree has the trivial semigroup
${\mathbb N}_{0}$ as its root and the children of each node are the semigroups
obtained by taking away one by one the right generators of the parent.
This construction was first considered in
\cite{RoGaGaJi:fundamentalgaps}.

Today, prior to this work, the RGD-algorithm \cite{rgd} is the fastest
implemented algorithm capable of computing $n_{g}$ without bounds on
$g$, the seeds algorithm \cite{seeds1,seeds2} is the fastest implemented
algorithm, but limited to genera up to around one half of the maximum size
of native integers. On its side, the Fromentin-Hivert algorithm
\cite{FromentinHivert} is the one that has been used to compute
$n_{g}$ up to a largest genus (up to genus 75)
\cite{DelgadoEliahouFromentin}.

Given a numerical semigroup $\Lambda $, consider the \emph{conductor}
$c(\Lambda )$, which is the first non-gap larger than the Frobenius number,
the \emph{multiplicity} $m(\Lambda )$, which is its first non-zero non-gap,
the \emph{jump} $u(\Lambda )$, which is the difference between the second
and first non-zero non-gaps, and the \emph{efficacy} $r(\Lambda )$, which
is the number of right generators. An \emph{encoding} of a numerical semigroup
will be a finite set of finite parameters which, together with its conductor,
multiplicity, jump, efficacy, and genus, uniquely define the numerical
semigroup. We will describe a general framework for algorithms exploring
the semigroup tree up to a given genus for a general encoding. This framework
fits the previous RGD-algorithm and the seeds algorithm. It just involves
the notions of right generators, strong generators (defined in Section~\ref{s:framework}),
pseudo-ordinary and quasi-ordinary semigroups.

Then we will present an encoding system that, used in the previous framework,
will allow to visit all nodes at a given depth while trimming all the leaves
and branches of ${\mathscr T}$ not arriving to that depth, providing the
so-called \emph{unleaved tree} of numerical semigroups of the given depth
(see \reftext{Fig.~\ref{fig:tresarbresunleaved}}). This represents a new paradigm
in the exploration of the semigroup tree up to a given genus, which will
lead to much more efficient algorithms than the ones used up to date. The
idea of pruning branches of the tree has been used in the context of the
search for counterexamples of Wilf's conjecture (or to related Eliahou's
conditions) in \cite{DelgadoTrimming,DelgadoEliahouFromentin}, where nodes
whose descendants are known to verify the conjecture are not further explored.

Our algorithm has already been used to compute the number of semigroups
of genus $76$ and $77$:
\begin{equation*}
n_{76}=29028294421710227
\end{equation*}
\begin{equation*}
n_{77}=47008818196495180
\end{equation*}

\section{General algorithms for visiting semigroups up to a given genus}
\label{s:framework}

Suppose that we can associate a unique finite set of finite parameters
$E(\Lambda )$ with each numerical semigroup $\Lambda $ such that the set
$E(\Lambda )$ together with $g(\Lambda )$, $c(\Lambda )$,
$m(\Lambda )$, $u(\Lambda )$, $r(\Lambda )$, completely determines
$\Lambda $. We call it an \emph{encoding} of $\Lambda $. Examples of encodings
are
\begin{itemize}
\item The minimal set of generators;
\item The \emph{Ap\'{e}ry set}, defined as the minimal elements of
$\Lambda $ of each congruence class modulo $m$, or the set of \emph{Kunz
coordinates}, defined as the integer quotients of the elements in the Ar\'{e}ry
set when divided by $m$;
\item The set of left elements or the set of gaps;
\item The sequence
$\nu _{i}=\#\{\lambda _{j}\in \Lambda :\lambda _{i}-\lambda _{j}\in
\Lambda \}$ up to $i=2c-g-1$ \cite{KiPe,Br:Acute};
\item The operation $\oplus $ defined by $i\oplus j=k$ if
$\lambda _{i}+\lambda _{j}=\lambda _{k}$ for all integers $i,j$ with
$0\leqslant i,j\leqslant 2c-g-1$ (indeed,
$\nu _{i}=\#\{(j,k)\in{\mathbb N}^{2}:j\oplus k=i\}$) \cite{Br:ANote};
\item The \emph{decomposition numbers}
$d_{i}=\lceil \frac{\nu _{i}}{2}\rceil $ up to $i=2c-g-1$
\cite{FromentinHivert};
\item The \emph{right-generators descendant} (RGD) of $\Lambda $, defined
as the numerical semigroup obtained by removing from $\Lambda $ all its
right generators, up to its $(c+m)$th element \cite{rgd};
\item The bitstreams encoding the gaps and the set of \emph{seeds} of order
up to $c-g$, where a seed of order $p$ is an element larger than the Frobenius
number and necessarily smaller than $c+\lambda _{p}-\lambda _{p-1}$, such
that $\lambda _{s}+\lambda _{p}\neq \lambda _{i}+\lambda _{j}$ for all
$p < i\leqslant j < s$ \cite{seeds1,seeds2}.
\end{itemize}

One can prove that all right generators of a semigroup $\Lambda $ with
conductor $c$ and multiplicity $m$ are at most $c+m-1$. Furthermore, it
was already proved in \cite{Br:bounds} that if
$\sigma _{1},\dots ,\sigma _{r}$ are the right generators of
$\Lambda $, with $\sigma _{1}<\sigma _{2}<\dots <\sigma _{r}$, then the
right generators of $\Lambda \setminus \{\sigma _{i}\}$ are either
$\{\sigma _{i+1},\dots ,\sigma _{r}\}$ or
$\{\sigma _{i+1},\dots ,\sigma _{r},\sigma _{i}+m\}$. We say that
$\sigma _{i}$ is a \emph{strong generator} of $\Lambda $ in the second case.
Otherwise, we say that it is a \emph{weak} generator. We will call
$\sigma _{1}$ the \emph{primogenial right generator} of $\Lambda $ and
$\Lambda \setminus \{\sigma _{1}\}$ the \emph{primogenial child} of
$\Lambda $. If $i\neq 1$, we say that $\sigma _{i-1}$ is the \emph{preceding
sibling} of $\sigma _{i}$ and $\Lambda \setminus \{\sigma _{i-1}\}$ is
the \emph{preceding sibling} of $\Lambda \setminus \{\sigma _{i}\}$.

Hence, for the exploration of ${\mathscr T}$ one is interested in encodings
of semigroups for which: there exists an efficient procedure to identify
right generators and strong right generators; there exists an efficient
procedure to obtain $E(\Lambda \setminus \{\sigma _{i}\})$, either from
the encoding of the parent $E(\Lambda )$ or from the encoding of the preceding
sibling $E(\Lambda \setminus \{\sigma _{i-1}\})$, should it exist (i.e.
if $i\neq 1$).

Furthermore, we have the following results from \cite{rgd}.
%
\begin{lemma}%
\label{l:threshold}
\cite[Lemma 3.1 and Lemma 3.3]{rgd} %
Let $\Lambda $ be a non-ordinary semigroup with multiplicity $m$ and jump
$u$. If $\sigma \geqslant c+u$, then $\sigma $ is not a strong generator.

Let $\Lambda $ be a pseudo-ordinary semigroup with multiplicity $m$ and
jump $u$. Then $c=m+u$ and all integers between $c$ and $c+m-1$ are right
generators except the integer $2m$. Furthermore, a right generator
$\sigma $ is strong if and only if $\sigma <c+u$.
\end{lemma}

Suppose that we have the following procedures
\begin{itemize}
\item \Call{CheckRightGenerator}{$E(\Lambda ),g,c,m,u,r,\sigma $}$=$
\Call{CRG}{$E(\Lambda ),g,c,m,u,r,\sigma $}, to check whether
$\sigma $ is a right generator of $\Lambda $;
\item \Call{CheckStrongGenerator}{$E(\Lambda ),g,c,m,u,r,\sigma $}$=$
\Call{CSG}{$E(\Lambda ),g,c,m,u,r,\sigma $}, to check whether a right generator
$\sigma $ is a strong generator of $\Lambda $;
\item \Call{EncodingFromParent}{$E(\Lambda ),g,c,m,u,r,\sigma $}$=$
\Call{EFP}{$E(\Lambda ),g,c,m,u,r,\sigma $}, to obtain
$E(\Lambda \setminus \{\sigma \})$, where $\sigma $ is a right generator
of $\Lambda $;
\item For $i\neq 1$,
\Call{EncodingFromPrecedingSibling}{$E(\Lambda \setminus \{\sigma _{i-1}\}),g,c,m,u,r,\sigma _{i}$}$=$
\Call{EFPS}{$E(\Lambda \setminus \{\sigma _{i-1}\}),g,\allowbreak c,m,u,r,\sigma _{i}$},
to obtain $E(\Lambda \setminus \{\sigma _{i}\})$, where
$\sigma _{i}$ is a non-primogenial right generator and
$\sigma _{i-1}$ is its preceding sibling.
\end{itemize}

Using these procedures and \reftext{Lemma~\ref{l:threshold}}, we will define a general
algorithm to visit semigroups of a given genus. We first present two subprocedures
and at the end the general procedure.

The procedure \Call{Descend}{} detailed in \reftext{Algorithm~\ref{a:descend}} recursively
explores the descendants of a non-ordinary and non-pseudo-ordinary semigroup
up to genus $\gamma $. Its input parameters are $E(\Lambda )$,
$m=m(\Lambda )$, $u=u(\Lambda )$, $c=c(\Lambda )$, $g=g(\Lambda )$,
$r=r(\Lambda )$, and $\gamma $.

\begin{algorithm}\caption{{\sc Descend}}\label{a:descend}
  \begin{algorithmic}[H]
    \Procedure{Descend}{$E(\Lambda),g,c,m,u,r,\gamma$}
    \State Visit $\Lambda$
    \If{$g<\gamma$}
    \State{$\tilde r\gets r$}
    \State{$\tilde E\gets E$}
    \For {$\sigma$ \From $c$ \To $c+u-1$}\Comment{{\color{blue}where strong generators may occur}}
    \If{\Call{CRG}{$\tilde E,g,c,m,u,\tilde r,\sigma$} $=$ \True}
    \If{$\sigma$ is primogenial}
    \State{$\tilde E\gets$ \Call{EFP}{$\tilde E,g,c,m,u,\tilde r,\sigma$}}
    \Else
    \State{$\tilde E\gets $ \Call{EFPS}{$\tilde E,g,c,m,u,\tilde r,\sigma$}}
    \EndIf
    \If{\Call{CSG}{$\tilde E,g,c,m,u,r,\sigma$} $=$ \True}
    \State {\Call{Descend}{$\tilde E,g+1,\sigma+1,m,u,\tilde r,\gamma$}}
    \State {$\tilde r\gets \tilde r-1$}
    \Else
    \State {$\tilde r\gets \tilde r-1$}
    \State {\Call{Descend}{$\tilde E,g+1,\sigma+1,m,u,\tilde r,\gamma$}}
    \EndIf
    \EndIf
    \EndFor
    \While{$\tilde r>0$}
    \Comment{{\color{blue}no more strong generators}}
    \If{\Call{CRG}{$\tilde E,g,c,m,u,\tilde r,\sigma$} $=$ \True}
    \If{$\sigma$ is primogenial}
    \State{$\tilde E\gets$ \Call{EFP}{$\tilde E,g,c,m,u,r,\sigma$}}
    \Else
    \State{$\tilde E\gets$ \Call{EFPS}{$\tilde E,g,c,m,u,r,\sigma$}}
    \EndIf
    \State $\tilde r\gets \tilde r-1$
    \State {\Call{Descend}{$\tilde E,g+1,\sigma+1,m,u,\tilde r,\gamma$}}
    \EndIf
    \EndWhile
    \EndIf
    \EndProcedure
  \end{algorithmic}
  \end{algorithm}

The procedure $\Call{PseudoDescend}{}$ explores the children of a given
pseudo-ordinary semigroup~$\Lambda $ that are not pseudo-ordinary, and
calls $\Call{Descend}{}$ to explore all their respective descendants in
${\mathscr T}$. It is detailed in \reftext{Algorithm~\ref{a:pseudo}}. In particular,
the variable $\sigma $ is set to never be $c$ since, otherwise,
$\Lambda \setminus \{\sigma \}$ would be another pseudo-ordinary semigroup.

\mut{
\begin{algorithm}%
\caption{\textsc{PseudoDescend}.}%
\label{a:pseudo}
\begin{algorithmic}
\Procedure{PseudoDescend}{$E(\Lambda ),c,m,u,r,\gamma $}
    \State Visit $\Lambda $
    \State{$\tilde r\gets r$}
    \State{$\tilde E\gets E$}
    \For {$\sigma $ \textbf{ from }$c+1$ \textbf{ to }$c+u-1$}\Comment{{\color{blue}where right generators are strong}}
    \If{$\sigma \neq 2m$}
    \State {\Call{Descend}{$\tilde E,c-1,\sigma +1,m,u,\tilde r,\gamma $}}
    \State $\tilde r\gets \tilde r-1$
    \EndIf
    \EndFor
    \For {$\sigma $ \textbf{ from }$c+u$ \textbf{ to }$c+m-1$}\Comment{{\color{blue}where right generators are weak}}
    \If{$\sigma \neq 2m$}
    \State $\tilde r\gets \tilde r-1$
    \State {\Call{Descend}{$\tilde E,c-1,\sigma +1,m,u,\tilde r,\gamma $}}
    \EndIf
    \EndFor
    \EndProcedure
  \end{algorithmic}%
\end{algorithm}
}

\begin{algorithm}\caption{{\sc PseudoDescend}}\label{a:pseudo}
  \begin{algorithmic}
    \Procedure{PseudoDescend}{$E(\Lambda),c,m,u,r,\gamma$}
    \State Visit $\Lambda$
    \State{$\tilde r\gets r$}
    \State{$\tilde E\gets E$}
    \For {$\sigma$ \From $c+1$ \To $c+u-1$}\Comment{{\color{blue}where right generators are strong}}
    \If{$\sigma\neq 2m$}
    \State {\Call{Descend}{$\tilde E,c-1,\sigma+1,m,u,\tilde r,\gamma$}}
    \State $\tilde r\gets \tilde r-1$
    \EndIf
    \EndFor
    \For {$\sigma$ \From $c+u$ \To $c+m-1$}\Comment{{\color{blue}where right generators are weak}}
    \If{$\sigma\neq 2m$}
    \State $\tilde r\gets \tilde r-1$
    \State {\Call{Descend}{$\tilde E,c-1,\sigma+1,m,u,\tilde r,\gamma$}}
    \EndIf
    \EndFor
    \EndProcedure
  \end{algorithmic}\end{algorithm}

Now we are ready to define a general exploring algorithm. For an integer
$m$ let $O_{m}$ be the unique ordinary semigroup of multiplicity $m$. For
$2\leqslant u\leqslant m$, let $P_{m,u}$ be the unique pseudo-ordinary
semigroup with multiplicity $m$ and jump $u$. For
$m+1\leqslant F\leqslant 2m-1$, let $Q_{m,F}$ be the unique quasi-ordinary
semigroup with multiplicity $m$ and Frobenius number $F$. Let
$H_{g}$ be the hyperelliptic semigroup of genus $g$, that is, the semigroup
generated by $2$ and $2g+1$.

Let ${\mathscr T}_{m}$ be the subtree of ${\mathscr T}$ with all the semigroups
of multiplicity $m$. It is the subtree that contains $O_{m}$ together with
all its descendants, except for the branch emerging from its unique ordinary
child, $O_{m+1}$.

In turn, ${\mathscr T}_{m}$ \xch{can be split}{can be splitted} into its $m-1$ subtrees
${\mathscr T}_{m,u}$, for $2\leqslant u\leqslant m$, each of which contains
all the semigroups of multiplicity $m$ and jump $u$ and the tree
${\mathscr Q}_{m}$ which is rooted in $O_{m}$ and contains all the quasi-ordinary
semigroups of multiplicity $m$ and all their respective descendants. That
is, ${\mathscr Q}_{m}$ contains all the semigroups with multiplicity
$m$ and jump $u=1$.

Notice that ${\mathscr T}_{m,u}$ is the subtree of
${\mathscr T}_{m}$ that contains $P_{m,u}$ together with all its descendants,
except for the branch emerging from its unique pseudo-ordinary child,
$P_{m,u+1}$, should $u<m$.

Combining the previous procedures \Call{Descend}{} and
\Call{PseudoDescend}{} we can define the algorithm
\Call{ExploreTree}{} detailed in \reftext{Algorithm~\ref{a:exploretree}} for exploring
${\mathscr T}$ up to a given genus $\gamma $. It separately explores the
trees ${\mathscr T}_{1},\dots ,{\mathscr T}_{\gamma +1}$, and, within the
exploration of each ${\mathscr T}_{m}$, it separately explores the tree
${\mathscr Q}_{m}$ and the trees ${\mathscr T}_{m,2}$,
${\mathscr T}_{m,3}$, \dots, ${\mathscr T}_{m,min}$, where
$min=\min \{m,\gamma +2-m\}$.

Hence, \Call{ExploreTree}{} can be parallelized by the multiplicity and
the jump in a \xch{straightforward}{straighforward} way. In turn, the exploration of
${\mathscr Q}_{m}$ can be parallelized by the $m$'th gap, that is, the
Frobenius number.

\mut{
\begin{algorithm}%
\caption{\textsc{ExploreTree}.}%
\label{a:exploretree}
\begin{algorithmic}
\Procedure{ExploreTree}{$\gamma $}
  \State {Visit ${\mathbb N}_{0}$} \Comment{{\color{blue} the unique semigroup with $m=1$}}
  \For{$g$ \textbf{ from }$1$ \textbf{ to }$\gamma $} \Comment{{\color{blue} the semigroups with $m=2$}}
  \State{Visit $H_{g}$}
  \EndFor
  \For {$m$ \textbf{ from }$3$ \textbf{ to }$\gamma $}
  \State{Visit $O_{m}$}
  \State{$min\gets \min \{m,\gamma +2-m\}$}
 \For{$u$ \textbf{ from }2 \textbf{ to }$min-1$}
 \State{Visit $P_{m,u}$}
 \State {\Call{PseudoDescend}{$E(P_{m,u}),m+u,m,u,m-2,\gamma $}} \Comment{{\color{blue}$c=m+u$, $r=m-2$}}
 \EndFor
  \State{Visit $P_{m,min}$}
  \If{$min<\gamma +2-m$}
  \State \Call{PseudoDescend}{$E(P_{m,m}),2m,m,m,m-1,\gamma $}  \Comment{{\color{blue}$u=m$, 
      $r=m-1$}}
  \EndIf
  \State{$r\gets m-3$}
  \For{$\sigma $ \textbf{ from }$m+2$ \textbf{ to }$2m-2$} \Comment{{\color{blue}loop on the quasi-ordinaries}}
 \State{\Call{Descend}{$E(Q_{m,\sigma}),m,\sigma +1,m,1,r,\gamma $}}
 \State{$r\gets r-1$}
 \EndFor
  \State{Visit $Q_{m,2m-1}$} \Comment{{\color{blue} $Q_{m,2m-1}$ has no descendants if $m>2$}}
  \EndFor
  \State{Visit $O_{\gamma +1}$} \Comment{{\color{blue} the unique semigroup with $g\leqslant \gamma $ and $m=\gamma +1$}}
  \EndProcedure
  \end{algorithmic}
\end{algorithm}
}

\begin{algorithm}\caption{{\sc ExploreTree}}\label{a:exploretree}
    \begin{algorithmic}
  \Procedure{ExploreTree}{$\gamma$}
  \State {Visit ${\mathbb N}_0$} \Comment{{\color{blue} the unique semigroup with $m=1$}}
  \For{$g$ \From $1$ \To $\gamma$} \Comment{{\color{blue} the semigroups with $m=2$}}
  \State{Visit $H_g$}
  \EndFor
  \For {$m$ \From $3$ \To $\gamma$}
  \State{Visit $O_m$}
  \State{$min\gets\min\{m,\gamma+2-m\}$}
 \For{$u$ \From 2 \To $min-1$}
 \State{Visit $P_{m,u}$}
 \State {\Call{PseudoDescend}{$E(P_{m,u}),m+u,m,u,m-2,\gamma$}} \Comment{{\color{blue}$c=m+u$, $r=m-2$}}
 \EndFor
  \State{Visit $P_{m,min}$}
  \If{$min<\gamma+2-m$}
  \State \Call{PseudoDescend}{$E(P_{m,m}),2m,m,m,m-1,\gamma$}  \Comment{{\color{blue}$u=m$, 
      $r=m-1$}}
  \EndIf
  \State{$r\gets m-3$}
  \For{$\sigma$ \From $m+2$ \To $2m-2$} \Comment{{\color{blue}loop on the quasi-ordinaries}}
 \State{\Call{Descend}{$E(Q_{m,\sigma}),m,\sigma+1,m,1,r,\gamma$}}
 \State{$r\gets r-1$}
 \EndFor
  \State{Visit $Q_{m,2m-1}$} \Comment{{\color{blue} $Q_{m,2m-1}$ has no descendants if $m>2$}}
  \EndFor
  \State{Visit $O_{\gamma+1}$} \Comment{{\color{blue} the unique semigroup with $g\leq \gamma$ and $m=\gamma+1$}}
  \EndProcedure
  \end{algorithmic}
\end{algorithm}

\section{Encoding by shrinking}
\label{sec3}

\subsection{A new encoding}
\label{sec3.1}

Given a numerical semigroup $\Lambda $ denote $L(\Lambda )$ its set of
\emph{left elements}, that is, its elements that are smaller than its Frobenius
number. Now, given a numerical semigroup $\Lambda $, define
\begin{eqnarray*}
\omega (\Lambda )&:=&\gcd (L(\Lambda ))
\\
\overline{\Lambda }&:=&\left \langle
\frac{L(\Lambda )}{\omega (\Lambda )}\right \rangle
\end{eqnarray*}

The tuple $\omega (\Lambda ),\overline{\Lambda }$ is an encoding of
$\Lambda $ as defined in Section~\ref{s:framework}. Indeed, together with
$c(\Lambda )$ it uniquely determines $\Lambda $, since
\begin{equation*}
\Lambda =\omega (\Lambda )\overline{\Lambda }\cup \left (c(\Lambda )+{
\mathbb N}_{0}\right ).
\end{equation*}

\begin{example}
Consider the semigroup
$\Lambda =\{0,8,10,16,17,\dots \}$. In this case,
$L(\Lambda )=\{0,8,10\}$. Hence, their $\gcd $ is $2$ and, so,
$\omega (\Lambda )=2$. Now,
$\frac{L(\Lambda )}{\omega (\Lambda )}=\{0,4,5\}$ and
$\overline{\Lambda }=\langle 4,5\rangle =\{0,4,5,8,9,10,12,\allowbreak 13,\dots
\}$.

Now,
$\omega (\Lambda )\overline{\Lambda }=\{0,8,10,16,18,20,24,26,28,
\dots \}$, while
$c(\Lambda )+{\mathbb N}_{0}=16+{\mathbb N}_{0}=\{16,17,\dots \}$. Hence,
the union of both sets is
\begin{equation*}
\{0,8,10,16,18,20,24,26,28,\dots \}\cup \{16,17,\dots \}=\{0,8,10,16,17,
\dots \}=\Lambda
\end{equation*}
\end{example}

In the definition of an encoding system we required its parameters to be
finite. Here we assume that $\overline{\Lambda }$ is represented by a finite
set, for instance, by its left elements. We call
$\overline{\Lambda }$ the \emph{shrinking} of $\Lambda $.

\subsection{\textsc{CheckRightGenerator} and \textsc{CheckStrongGenerator}}
\label{sec3.2}

Next we show that we have a direct way to determine from
$\overline{\Lambda }$ and $\omega (\Lambda )$ the right generators of
$\Lambda $, and that it is \xch{equally}{equaly} easy to identify the strong generators.

\begin{lemma}
\label{lem3.1}
Let $\Lambda $ be a numerical semigroup and let $c=c(\Lambda )$,
$m=m(\Lambda )$, $u=u(\Lambda )$, $\omega =\omega (\Lambda )$.
\begin{enumerate}
\item Let $c\leqslant \sigma <c+m$. The element
$\sigma \in \Lambda $ is a right generator of $\Lambda $ if and only if
either
\begin{enumerate}%
\item[(a)] $\sigma \not \equiv 0\mod \omega $
\item[(b)] $\sigma \equiv 0\mod \omega $ and
$\frac{\sigma}{\omega}\not \in \overline{\Lambda }$.
\end{enumerate}
\item In case (a), $\sigma $ is strong if and only if $\sigma <c+u$.
\item In case (b), $\sigma $ is strong if and only if
$\frac{\sigma +m}{\omega (\Lambda \setminus \{\sigma \})}\not \in
\overline{\Lambda \setminus \{\sigma \}}$.
\end{enumerate}
\end{lemma}

\begin{proof}
It is obvious that if $\sigma \not \equiv 0\mod \omega $ then
$\sigma $ is a right generator. If $\sigma \equiv 0\mod \omega $, then
it is not a right generator if and only if it is generated by the left
elements, which is equivalent to
$\frac{\sigma}{\omega}\in \overline{\Lambda }$.

For the second statement, on the one hand it follows from \reftext{Lemma~\ref{l:threshold}}
that if $\sigma \geqslant c+u$ then $\sigma $ is not strong. Now, suppose
that $\sigma <c+u$ and suppose that $\sigma $ is not strong, that is,
$\sigma +m=a+b$ with $a,b\in\Lambda$, $\{a,b\}\neq \{\sigma ,m\}$, and
$\{a,b\}\neq \{0,\sigma +m\}$. Since $a,b\neq 0,m$, then
$a,b\geqslant m+u>m+\sigma -c$ and, by the equality $a+b=\sigma +m$ we
deduce
$\left \{
\begin{array}{l}
\sigma +m>a+m+\sigma -c
\\
\sigma +m>b+m+\sigma -c%
\end{array}
\right .$ and, so, $a,b<c$. Then, $a,b\equiv 0\mod \omega $ and
$a+b\equiv 0\mod \omega $ implying that $\sigma \equiv 0\mod\omega $, a
contradiction.

For the third item, we know that $\sigma $ is strong if and only if
$\sigma +m$ is a minimal generator of
$\Lambda \setminus \{\sigma \}$, and, by the first item, this is equivalent
to
$\frac{\sigma +m}{\omega (\Lambda \setminus \{\sigma \})}\not \in
\overline{\Lambda \setminus \{\sigma \}}$.
\end{proof}

\subsection{Encoding pseudo-ordinary and quasi-ordinary semigroups and \textsc{EncodingFromParent} and \textsc{EncodingFromPrecedingSibling}}
\label{sec3.3}

\paragraph*{Semigroups generated by an interval}
In order to explain the encoding procedures, we need some results on semigroups
generated by intervals. Let $\Lambda _{\{i,\dots ,j\}}$ be the semigroup
generated by the interval $\{i,i+1,\dots ,j\}$. Notice that it is the union
of sets $S_{k}=\{ki,\dots ,kj\}$.

\begin{lemma}
\label{lem3.2}
The conductor of $\Lambda _{\{i,\dots ,j\}}$ is
$i\lfloor \frac{j-2}{j-i}\rfloor $. The genus of
$\Lambda _{\{i,\dots ,j\}}$ is
$\sum _{k=1}^{i\lfloor \frac{j-2}{j-i}\rfloor}(i+(k-1)(i-j)-1)$.
\end{lemma}

\begin{proof}
Let $D_{1}$ be the set of gaps between $0$ and $S_{1}$ and, in general,
$D_{k}$ the set of gaps between $S_{k-1}$ and $S_{k}$. One can check that
$\#D_{k}=i+(k-1)(i-j)-1$ as long as this is a positive amount, although
we will also extend this equality for the case of negative cardinality.
Let $k_{c}$ be such that the conductor of $\Lambda $ is the smallest element
in $S_{k_{c}}$. Equivalently,
%
\begin{equation}
\label{e:kc}
k_{c}=\frac{c}{i}.%
\end{equation}

Notice that the Frobenius number of $\Lambda _{\{i,\dots ,j\}}$ is the
largest element of $D_{k_{c}}$ and $D_{i}\neq 0$ if and only if
$i\leqslant k_{c}$. Hence, $k_{c}$ is the largest element such that
$\#D_{k}\geqslant 1$, that is,
\begin{eqnarray*}
i+(k_{c}-1)(i-j)-1&\geqslant &1
\\
i+(k_{c}-1)(i-j)&\geqslant &2
\\
k_{c}&\leqslant &\frac{2-i}{i-j}+1=\frac{2-j}{i-j}=\frac{j-2}{j-i}
\end{eqnarray*}
So,
\begin{equation*}
k_{c}=\left \lfloor \frac{j-2}{j-i}\right \rfloor .
\end{equation*}
Now the result follows from \reftext{\eqref{e:kc}}.
\end{proof}

\paragraph*{Encoding pseudo-ordinary and quasi-ordinary semigroups}
It is easy to check that $\omega (P_{m,u})=m$ and
$\overline{P}_{m,u}={\mathbb N}_{0}$. In particular
$Q_{m,m+1}=P_{m,2}$ and so $\omega (Q_{m,m+1})=m$ and
$\overline{Q}_{m,m+1}={\mathbb N}_{0}$. Similarly, if $F>m+1$, then
$L(Q_{m,F})=\{m,m+1,\dots ,F-1\}$ and, so $\omega (Q_{m,F})=1$ and
$\overline{Q}_{m,F}=\Lambda _{\{m,\dots ,F-1\}}$.

\paragraph*{Encoding from the parent}
Let us see how we can encode a semigroup from the encoding of its parent.
If $A,B$ are two sets, we will use $A+B$ to refer to the set
$\{a+b:a\in A,b\in B\}$. If $\Lambda _{1},\Lambda _{2}$ are submonoids
of ${\mathbb N}_{0}$, then $\Lambda _{1}+\Lambda _{2}$ is also a submonoid
and, if either $\Lambda _{1}$ or $\Lambda _{2}$ is a numerical semigroup,
then $\Lambda _{1}+\Lambda _{2}$ is a numerical semigroup and, indeed,
it is the minimum semigroup that contains $\Lambda _{1}$ and
$\Lambda _{2}$.

\begin{lemma}%
\label{l:efp}
Let $\Lambda $ be a numerical semigroup with conductor $c$. Let
$\omega =\omega (\Lambda )$,
$\tilde\Lambda =\Lambda \setminus \{\sigma \}$ and
$\tilde\omega =\omega (\tilde\Lambda )$.
\begin{itemize}
\item If $\sigma =c$,
\begin{itemize}
\item $\tilde\omega =\omega $,
\item $\overline{\tilde \Lambda }=\overline{\Lambda }$,
\item $c(\overline{\tilde \Lambda })=c(\overline{\Lambda })$.
\end{itemize}
\item If $\sigma =c+1$,
\begin{itemize}
\item $\tilde\omega =\gcd (\omega ,c)$,
\item
$\overline{\tilde \Lambda }=\frac{\omega}{\tilde\omega}
\overline{\Lambda }+\frac{c}{\tilde\omega}{\mathbb N}_{0}$,
\item
$c(\overline{\tilde \Lambda })\leqslant c(\overline{\Lambda })
\frac{\omega}{\tilde\omega}+(\frac{c}{\tilde\omega}-1)(
\frac{\omega}{\tilde\omega}-1)$.
\end{itemize}
\item If $\sigma >c+1$,
\begin{itemize}
\item $\tilde\omega =1$,
\item
$\overline{\tilde \Lambda }=\omega \overline{\Lambda }+\Lambda _{c,
\dots ,{\sigma -1}}$,
\item
$c(\overline{\tilde \Lambda })\leqslant c(\overline{\Lambda })\omega +(
\lfloor \frac{\omega -2}{\sigma -c-1}\rfloor +1)c$
\end{itemize}
\end{itemize}
\end{lemma}

\begin{proof}
The case $\sigma =c$ is clear.

Suppose $\sigma =c+1$. In this case $\overline{\tilde \Lambda }$ contains
$\frac{\omega}{\tilde\omega}\overline\Lambda$ and
$\frac{c}{\tilde\omega}$. Consider the list
\begin{equation*}
c(\overline{\Lambda })\frac{\omega}{\tilde\omega},\ \ \ c(
\overline{\Lambda })\frac{\omega}{\tilde\omega}+
\frac{c}{\tilde\omega},\ \ \ c(\overline{\Lambda })
\frac{\omega}{\tilde\omega}+2\frac{c}{\tilde\omega},\ \ \ c(
\overline{\Lambda })\frac{\omega}{\tilde\omega}+3
\frac{c}{\tilde\omega},\ \ \ \dots \ \ \ c(\overline{\Lambda })
\frac{\omega}{\tilde\omega}+\left (\frac{\omega}{\tilde\omega}-1
\right )\frac{c}{\tilde\omega}.
\end{equation*}
It contains $\frac{\omega}{\tilde\omega}$ elements of
$\overline{\tilde \Lambda }$ that are not congruent modulo
$\frac{\omega}{\tilde\omega}$ because $\frac{c}{\tilde\omega}$ and $\frac{\omega}{\tilde\omega}$ are coprime. Hence, the elements in the list range all
congruence classes modulo $\frac{\omega}{\tilde\omega}$. But, since
$(c(\overline{\Lambda })+i)\frac{\omega}{\tilde\omega}$ also belongs to
$\overline{\tilde \Lambda }$ for any positive integer $i$, the list
\begin{equation*}
(c(\overline{\Lambda })+i)\frac{\omega}{\tilde\omega},\ \ \ (c(
\overline{\Lambda })+i)\frac{\omega}{\tilde\omega}+
\frac{c}{\tilde\omega},\ \ \ (c(\overline{\Lambda })+i)
\frac{\omega}{\tilde\omega}+2\frac{c}{\tilde\omega},\ \ \ %
\ \ \ \dots \ \ \ (c(\overline{\Lambda })+i)
\frac{\omega}{\tilde\omega}+\left (\frac{\omega}{\tilde\omega}-1
\right )\frac{c}{\tilde\omega}
\end{equation*}
is contained in $\overline{\tilde \Lambda }$. Hence, any element in the first list plus any multiple of $\frac{\omega}{\tilde\omega}$ is contained in $\overline{\tilde \Lambda }$. Now, the elements between
$c(\overline{\Lambda })\frac{\omega}{\tilde\omega}+(
\frac{\omega}{\tilde\omega}-1)\frac{c}{\tilde\omega}-
\frac{\omega}{\tilde\omega}+1$ and
$c(\overline{\Lambda })\frac{\omega}{\tilde\omega}+(
\frac{\omega}{\tilde\omega}-1)\frac{c}{\tilde\omega}$ all belong to
$\overline{\tilde \Lambda }$ and, so
$c(\overline{\tilde \Lambda })\leqslant c(\overline{\Lambda })
\frac{\omega}{\tilde\omega}+(\frac{\omega}{\tilde\omega}-1)
\frac{c}{\tilde\omega}-\frac{\omega}{\tilde\omega}+1=c(
\overline{\Lambda })\frac{\omega}{\tilde\omega}+(
\frac{\omega}{\tilde\omega}-1)(\frac{c}{\tilde\omega}-1)$.

Suppose $\sigma >c+1$. It is obvious that $\tilde\omega =1$. The semigroup
$\overline{\tilde \Lambda }$ must contain
$\Lambda _{c,\dots ,\sigma -1}$ and, so it contains the union of sets
$S_{k}=\{kc,\dots ,k(\sigma -1)\}$, which have $k(\sigma -1-c)+1$ consecutive
elements. In particular, if
$\ell =\lceil \frac{\omega -1}{\sigma -1-c}\rceil $, the set
$S_{\ell}$ contains $\omega $ consecutive elements. From this we deduce
that the intervals $(c(\overline{\Lambda })+i)\omega +S_{\ell}$ for all $i\geqslant 0$ contain
$\omega $ consecutive elements and they cover all the integers greater
than or equal to $c(\overline{\Lambda })\omega +\ell c$. Hence, the conductor
of $\overline{\tilde \Lambda }$ is at most
$c(\overline{\Lambda })\omega +\ell c$ and the result follows.
\end{proof}

\begin{example}
Consider the semigroup $\Lambda =\{0,12,16,18,19,\dots \}$. It satisfies
\begin{itemize}
\item $L(\Lambda )=\{0,12,16\}$,
\item $\omega =\gcd (0,12,16)=4$,
\item $\overline{\Lambda }=\langle 3,4\rangle =\{0,3,4,6,7\dots \}$,
\item $c(\overline{\Lambda })=6$.
\end{itemize}
Let us now see what happens when we remove $\sigma $ from $\Lambda $ in
the next cases:
\begin{itemize}
\item If $\sigma =c=18$, then
$\tilde\Lambda =\Lambda \setminus \{18\}= \{0,12,16,19,\dots \}$. In this
case,
\begin{itemize}
\item $L(\tilde\Lambda )=\{0,12,16\}=L(\Lambda )$,
\item
$\tilde\omega =\gcd (L(\tilde\Lambda ))=\gcd (L(\Lambda ))=\omega =4$,
\item
$\overline{\tilde \Lambda }=\langle
\frac{L(\tilde\Lambda )}{\tilde\omega}\rangle =\langle
\frac{L(\Lambda )}{\omega}\rangle =\overline{\Lambda }$.
\item Obviously, $c(\overline{\tilde\Lambda})=c(\overline\Lambda)$.
\end{itemize}
\item If $\sigma =c+1=19$, then
$\tilde\Lambda =\Lambda \setminus \{19\}= \{0,12,16,18,20,21,\dots \}$.
In this case,
\begin{itemize}
\item $L(\tilde\Lambda )=\{0,12,16,18\}$,
\item
$\tilde\omega =\gcd (L(\tilde\Lambda ))%
=\gcd (\gcd (0,12,16),18)=\gcd (\omega ,18)=2$,
\item
$\overline{\tilde \Lambda }=\left \langle
\frac{L(\tilde\Lambda )}{\tilde\omega}\right \rangle =\langle 0,6,8,9
\rangle =\{0,6,8, 9,%
12,%
14,%
15,%
16,%
17,%
18, %
20, 21, 22, 23,%
24,\dots \}$. Any even number in this set belongs to
$\frac{\omega}{\tilde\omega}\overline{\Lambda }=
\frac{\omega}{\tilde\omega}\overline{\Lambda }+\{0\}\subseteq
\frac{\omega}{\tilde\omega}\overline{\Lambda }+\frac{c}{\tilde\omega}{
\mathbb N}_{0}$, while any odd number belongs to
$\frac{\omega}{\tilde\omega}\overline{\Lambda }+\{9\}\subseteq
\frac{\omega}{\tilde\omega}\overline{\Lambda }+\frac{c}{\tilde\omega}{
\mathbb N}_{0}$, since $\frac{c}{\tilde\omega}=9$. Hence,
$\overline{\tilde \Lambda }=\frac{\omega}{\tilde\omega}
\overline{\Lambda }+\frac{c}{\tilde\omega}{\mathbb N}_{0}$.
\item $c(\overline{\tilde \Lambda })=20$ which coincides with the upper
bound in the lemma, that is, $6\cdot 2+(18/2-1)(2-1)=20$.
\end{itemize}
\item If $\sigma =c+2=20$, then
$\tilde\Lambda =\Lambda \setminus \{20\}=\{0,12,16,18,19,21,22,\dots
\}$. In this case,
\begin{itemize}
\item $L(\tilde\Lambda )=\{0,12,16,18,19\}$,
\item
$\tilde\omega =\gcd (L(\tilde\Lambda ))%
=\gcd (0,12,16,18,19)=1$,
\item
$\overline{\tilde \Lambda }=\langle L(\tilde\Lambda )\rangle =
\langle 0,12,16,18,19\rangle =\{0,{\mathbf{12}},{\mathbf{16}},18,19,%
{\mathbf{24}},%
{\mathbf{28}},%
30,%
31,%
{\mathbf{32}},%
34,%
35$,%

${\mathbf{36}},%
37,%
38,%
{\mathbf{40}},%
42,%
43,%
{\mathbf{44}},%
46,%
47,%
{\mathbf{48}},%
49,%
50,%
51,%
{\mathbf{52}},%
53,%
54,%
55,%
\dots \}$. We boldfaced the elements that belong to
$4\overline{\Lambda }=4\langle 3,4\rangle $. It can be easily checked that
the rest of the elements can be obtained as a sum of an element in
$4\langle 3,4\rangle $ plus a non-zero element in
$\langle 18,19\rangle $. Hence,
$\overline{\tilde \Lambda }=\omega \overline{\Lambda }+\Lambda _{18,19}$.
\item $c(\overline{\tilde \Lambda })=46$, while the upper bound is
$6\cdot 4+\left (\left \lfloor \frac{4-2}{20-18-1}\right \rfloor +1
\right )18=78$.%
\end{itemize}
\end{itemize}
\end{example}

\paragraph*{Encoding from the preceding sibling}
It is now straightforward to prove the following lemma.

\begin{lemma}%
\label{l:efps}%
Let $\Lambda $ be a numerical semigroup with conductor $c$. Suppose that
the right generators of $\Lambda $ are
$\sigma _{1}<\dots <\sigma _{r}$. Suppose that
$2\leqslant i\leqslant r$ and that $\sigma _{i-1}\neq c$. Let
$\tilde\Lambda =\Lambda \setminus \{\sigma _{i}\}$, and let
$\tilde\Lambda '=\Lambda \setminus \{\sigma _{i-1}\}$. If
$\omega (\tilde{\Lambda}')=1$, then
\begin{itemize}
\item[]
\begin{itemize}
\item $\tilde\omega =1$,
\item
$\overline{\tilde \Lambda }={{\overline{\tilde{\Lambda }}}'}+\Lambda _{
\sigma _{i-1},\dots ,{\sigma _{i}-1}}$,
\item
$c(\overline{\tilde \Lambda })\leqslant c({
\overline{\tilde{\Lambda }}}')$
\end{itemize}
\end{itemize}
\end{lemma}

Notice that condition $\omega (\tilde{\Lambda}')=1$ is satisfied whenever
$\sigma _{i-1}\geqslant c+2$.

\section{The unleaved tree}
\label{sec4}

The interest in encoding by the gcd and shrinking is the next result, which
was first proved in \cite{BrBu}, although in another context and with different
notation.

\begin{lemma}%
\label{lem4.1}
\cite[Theorem 10]{BrBu} If $\omega (\Lambda )\neq 1$, then
$\Lambda $ has descendants of any given genus. If
$\omega (\Lambda )=1$, the maximum genus of the descendants of
$\Lambda $ is the genus of $\overline{\Lambda }$ and there is only one
descendant of that genus.
\end{lemma}

Hence, if we only want to visit the semigroups of genus $\gamma $, we can
trim any numerical semigroup $\Lambda $ of $\mathscr{T}$ satisfying that
$\omega (\Lambda )=1$ and that the genus of $\overline{\Lambda }$ is smaller
than $\gamma $, together with all its \xch{descendants}{descedants}. We call \emph{unleaved
tree of genus $\gamma $} the subtree that we obtain. This tree has no leaves
of genus smaller than $\gamma $ and all its leaves are exactly the semigroups
of genus $\gamma $. In contrast, we will call \emph{complete tree of genus
$\gamma $} the subtree of ${\mathscr{T}}$ that contains all numerical semigroups
of genus up to~$\gamma $.

In \reftext{Fig.~\ref{fig:tresarbresunleaved}} we represent the complete tree of
genus $9$ (left), the unleaved tree of genus $9$ (center), and the edges
of the complete tree that do not appear in the unleaved tree (right).

\mut{
\begin{figure}
\caption{Complete tree, unleaved tree, and their difference for genus
$\gamma =9$.}%
\label{fig:tresarbresunleaved}
\end{figure}
}
  
\begin{figure}
  \caption{Complete tree, unleaved tree, and their difference for genus $\gamma=9$}\label{fig:tresarbresunleaved}
    \resizebox{\textwidth}{!}{\includegraphics{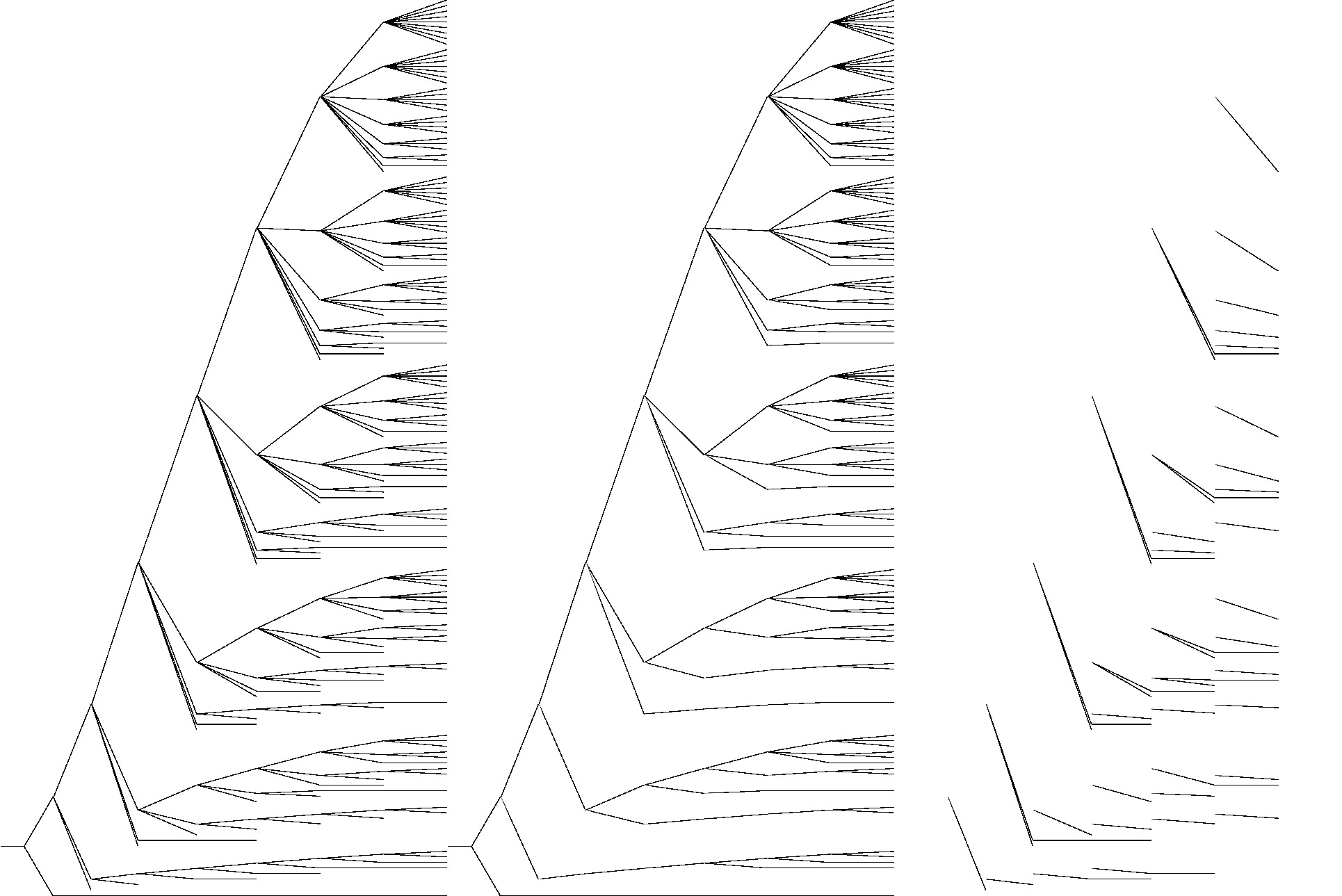}}
\end{figure}

If our objective in the exploration of the tree is counting, we can also
avoid visiting the nodes such that the genus of
$\overline{\Lambda }$ is exactly equal to $\gamma $ and just count~$1$.

A fact that makes trimming very efficient is that it is proved that most
numerical semigroups belong to finite chains \cite{rarity}, that is, most
numerical semigroups have the gcd of their left elements equal to
$1$.

Before deciding whether we trim a semigroup and all its descendants we
need to encode it, which takes some computing time. Thus, an important
trick of our algorithm is not encoding nodes that we know a priori that
will be trimmed. The result in the next lemma is in this direction.

\begin{lemma}
\label{lem4.2}
If $\Lambda $, $\Lambda '$ are siblings in the semigroup tree with
$F(\Lambda ')<F(\Lambda )$, then if $\Lambda '$ has no descendants of genus
$\gamma $ then neither does $\Lambda $.
\end{lemma}

\begin{proof}
If $\Lambda '$ has no descendants then $\omega (\Lambda ')=1$. The left
elements of $\Lambda $ are the left elements of $\Lambda '$ together with
the interval $\{F(\Lambda '),\dots ,F(\Lambda )-1\}$, which is not empty
and which contains elements not in $\Lambda '$. So,
$\omega (\Lambda )=\omega (\Lambda ')=1$ and $\overline{\Lambda '}$ is
strictly contained in $\overline{\Lambda }$. So,
$g(\overline{\Lambda '})>g(\overline{\Lambda })$. Hence, if
$g(\overline{\Lambda '})<\gamma $, then
$g(\overline{\Lambda })<\gamma $.
\end{proof}

All these results suggest the function \Call{DescendAndTrim}{} shown in
\reftext{Algorithm~\ref{a:descendandtrim}} in substitution of \Call{Descend}{}, when
the exploration of the tree is aimed at counting semigroups of genus
$\gamma $. If the exploration of the tree is aimed at visiting all nodes
of genus $\gamma $, and not just counting, then a statement for visiting
$\Lambda $ should be placed at the beginning, all the lines with the parameter
count, as well as all lines from line~\ref{line:ggran} to the end should
be omitted, and line~\ref{line:glimit} should be replaced by ``\textbf{if}
$g<\gamma $ \textbf{then}''.

Notice that the procedure \Call{PseudoDescend}{} can be equally optimized
using the same idea of trimming. We call
\Call{PseudoDescendAndTrim}{} the new function. For the sake of brevity
we do not put here the associated pseudo-code.

  \begin{algorithm}\caption{{\sc DescendAndTrim}}\label{a:descendandtrim}
  \begin{algorithmic}[1]{\tiny
    \Function{DescendAndTrim}{$\omega,\shrink{\Lambda},g,c,m,u,r,\gamma$}
    \State{count $\gets 0$} \label{line:count}
    \State{$\tilde r\gets r$}
    \State{$\tilde\omega,\shrink{\tilde\Lambda}\gets \omega,\shrink\Lambda$}
    \If{$g<\gamma-2$} \label{line:glimit}
    \State{keepgoing $\gets$ \True}
    \For {$\sigma$ \From $c$ \To $c+u-1$, {\bf while} keepgoing }
    \If{$\sigma\not\equiv 0\mod \omega$ \Or $\frac{\sigma}{\omega}\not\in\shrink\Lambda$ }
    \If{$\sigma$ is primogenial \Or the preceding sibling is $c$}
   \State{$\tilde\omega,\shrink{\tilde{\Lambda}}\gets$ \Call{EFP}{$\omega,\shrink{\Lambda},g,c,m,u,\tilde r,\sigma$}} \Comment{{\color{blue} Lemma~\ref{l:efp}}}
    \Else
   \State{$\tilde\omega,\shrink{\tilde\Lambda}\gets $ \Call{EFPS}{$\tilde\omega,\shrink{\tilde\Lambda},g,c,m,u,\tilde r,\sigma$}} \Comment{{\color{blue} Lemma~\ref{l:efps}}}
    \EndIf
    \If{$\tilde\omega\neq 1$}
    \If{$\sigma\not\equiv0\mod\tilde\omega$ \Or $\frac{\sigma+m}{\tilde\omega}\not\in\shrink{\tilde\Lambda}$}
    \State {count $\gets$ count$+$\Call{DescendAndTrim}{$\tilde \omega,\shrink{\tilde\Lambda},g+1,\sigma+1,m,u,\tilde r,\gamma$}}
    \State {$\tilde r\gets \tilde r-1$}
    \Else
    \State {$\tilde r\gets \tilde r-1$}
    \State {count $\gets$ count$+$\Call{DescendAndTrim}{$\tilde \omega,\shrink{\tilde\Lambda},g+1,\sigma+1,m,u,\tilde r,\gamma$}}
    \EndIf
    \Else
    \If{genus$(\shrink{\tilde{\Lambda}})\leq \gamma$}
    \If{genus$(\shrink{\tilde{\Lambda}})= \gamma$}
    \State{count $\gets$ count$+1$}
    \EndIf
    \State{keepgoing $\gets$ \False}
    \Else
    \If{$\sigma+m\not\in\shrink{\tilde\Lambda}$}
    \State {count $\gets$ count$+$\Call{DescendAndTrim}{$\tilde \omega,\shrink{\tilde\Lambda},g+1,\sigma+1,m,u,\tilde r,\gamma$}}
    \State {$\tilde r\gets \tilde r-1$}
    \Else
    \State {$\tilde r\gets \tilde r-1$}
    \State {count $\gets$ count$+$\Call{DescendAndTrim}{$\tilde \omega,\shrink{\tilde\Lambda},g+1,\sigma+1,m,u,\tilde r,\gamma$}}
    \EndIf\EndIf
    \EndIf
    \EndIf
    \EndFor
    \While{keepgoing {\textbf{ and }} $\tilde r>1$} \\\Comment{{\color{blue}no strong generators here, so if $\tilde r=1$ there will be no grand-children}}
    \If{$\sigma\not\equiv 0\mod \omega$ \Or $\frac{\sigma}{\omega}\not\in\shrink\Lambda$ }
    \If{$\sigma$ is primogenial \Or the preceding sibling is $c$}
    \State{$\tilde\omega,\tilde\Lambda\gets$ \Call{EFP}{$\tilde \omega,\shrink{\tilde\Lambda},g,c,m,u,r,\sigma$}} \Comment{{\color{blue} Lemma~\ref{l:efp}}}
    \Else
    \State{$\tilde\omega,\tilde\Lambda\gets$ \Call{EFPS}{$\tilde \omega,\shrink{\tilde\Lambda},g,c,m,u,r,\sigma$}} \Comment{{\color{blue} Lemma~\ref{l:efps}}}
    \EndIf
    \If{genus$(\shrink{\tilde{\Lambda}})\leq \gamma$}
    \If{genus$(\shrink{\tilde{\Lambda}})= \gamma$}
    \State{count $\gets$ count$+1$}
    \EndIf
    \State{keepgoing $\gets$ \False}
    \Else
    \State $\tilde r\gets \tilde r-1$
    \State {count $\gets$ count$+$\Call{DescendAndTrim}{$\tilde\omega,\shrink{\tilde\Lambda},g+1,\sigma+1,m,u,\tilde r,\gamma$}}
    \EndIf\EndIf
    \EndWhile
    \Else \label{line:ggran}\Comment{{\color{blue} if $g=\gamma-2$}}
    \For {$\sigma$ \From $c$ \To $c+u-1$}
    \If{$\sigma\not\equiv 0\mod \omega$ \Or $\frac{\sigma}{\omega}\not\in\shrink\Lambda$ }
    \If{$\sigma$ is primogenial \Or the preceding sibling is $c$}
    \State{$\tilde\omega,\tilde\Lambda\gets$ \Call{EFP}{$\tilde \omega,\shrink{\tilde\Lambda},g,c,m,u,r,\sigma$}} \Comment{{\color{blue} Lemma~\ref{l:efp}}}
    \Else
    \State{$\tilde\omega,\tilde\Lambda\gets$ \Call{EFPS}{$\tilde \omega,\shrink{\tilde\Lambda},g,c,m,u,r,\sigma$}} \Comment{{\color{blue} Lemma~\ref{l:efps}}}
    \EndIf
    \If{$\sigma\not\equiv0\mod\tilde\omega$ \Or $\frac{\sigma+m}{\tilde\omega}\not\in\shrink{\tilde\Lambda}$}
    \State{count $\gets$ count$+r$}
    \State {$\tilde r\gets \tilde r-1$}
    \Else
    \State {$\tilde r\gets \tilde r-1$}
    \State{count $\gets$ count$+r$}
    \EndIf\EndIf
    \EndFor
    \While{$\tilde r>1$}
    \State {$\tilde r\gets \tilde r-1$}
    \State{count $\gets$ count$+r$}
    \EndWhile
    \EndIf
    \State{\Return count}
    \EndFunction
}  \end{algorithmic}
\end{algorithm}

\section{More unvisited nodes}
\label{sec5}

\begin{lemma}%
\label{l:countgrandchildren}
\cite[Section 5]{rgd} %
The number of children of $P_{m,u}$ is $m-1$.

The number of grandchildren of $P_{m,u}$ is
\begin{equation*}
\left \{
\begin{array}{l@{\quad}l}
\binom{m-1}{2}+u &\mbox{ if }2u\leqslant m
\\
\binom{m-1}{2}+u -1&\mbox{ otherwise}
\end{array}
\right .
\end{equation*}
\end{lemma}

Notice that $P_{m,u'}$ is a descendant of $P_{m,u}$ as far as $u'>u$. If
we only want to count the semigroups of genus $\gamma $, in our algorithm
we need only to descend the pseudo-ordinary semigroups $P_{m,u}$ for
$u\leqslant \gamma -m-1$, as far as $\gamma -m\leqslant m$. Indeed,
$P_{m,\gamma -m+1}$ and $P_{m,\gamma -m+2}$ (should they exist) are descendants
of $P_{m,\gamma -m}$, and the descendants of $P_{m,\gamma -m}$ of genus
$\gamma $ are exactly its grandchildren, which can be counted by the formula
in \reftext{Lemma~\ref{l:countgrandchildren}}, without needing to visit them.

From another perspective, Rosales proved the result in the next lemma in
\cite{Rosales:m3}.

\begin{lemma}%
\label{l:minmultiplicities}
\cite[Corollary 10]{Rosales:m3} %
There is exactly one semigroup of genus $\gamma $ and multiplicity
$2$ and there are $\lfloor \frac{\gamma}{3}\rfloor +1$ semigroups of genus
$\gamma $ and multiplicity $3$
\end{lemma}

\begin{lemma}%
\label{l:maxmultiplicities}
If $\gamma \geqslant 8$, then there are
$\binom{\gamma -4}{4}+\binom{\gamma -2}{3}+\binom{\gamma -5}{2}+6
\gamma -14$ numerical semigroups of genus $\gamma $ and multiplicity larger
than or equal to $\gamma -3$.
\end{lemma}

\begin{proof}
It was proved in \cite[Theorem 9]{seeds2} that for $m\geqslant 4$, the
ordinary semigroup of multiplicity $m$ has exactly
$\binom{m}{3}+3m+3$ great-grandchildren. This implies that the ordinary
semigroup of multiplicity $m=\gamma -2$ has
$\binom{\gamma -2}{3}+3\gamma -3$ descendants of genus $\gamma $, and,
so, there are $\binom{\gamma -2}{3}+3\gamma -3$ semigroups of genus
$\gamma $ and multiplicity larger than or equal to $\gamma -2$.

Now we claim that there are \xch{exactly}{excatly}
$\binom{\gamma -4}{4}+\binom{\gamma -5}{2}+3\gamma -11$ semigroups of genus
$\gamma $ and multiplicity exactly equal to $\gamma -3$.

Suppose that a semigroup has genus $\gamma $ and multiplicity
$m=\gamma -3$. In particular, the semigroup has exactly four gaps larger
than the multiplicity. Define $A_{1}$ as the set of congruence classes
modulo $m$ of its gaps between $m$ and $2m$ and define subsequently
$A_{i}$ as the set of classes modulo $m$ of its gaps between $i\,m$ and
$(i+1)m$. It is obvious that $A_{i+1}\subseteq A_{i}$. Let $j$ be the maximum
integer such that $A_{j}\neq \emptyset $. The genus of the semigroup is
$m+3$ and the Frobenius number is $j\,m+\max (A_{j})$, from which
$j\,m+\max (A_{j})\leqslant 2m+5$ and, so
$1\leqslant \max (A_{j})\leqslant (2-j)m+5$. This implies that
$j\leqslant 2+\frac{4}{m}$ and so, for $m>4$ (which is the case if
$m=\gamma -3$ and $\gamma \geqslant 8$), we have $j\leqslant 2$. In particular,
$\#A_{1}+\#A_{2}=4$ and, if $A_{2}\neq \emptyset $, then
$\max (A_{2})\leqslant 5$.

Consider all the possible cases \xch{depending}{dependning} on the cardinality of
$A_{1}$. If $\#A_{1}=4$, then there are no restrictions on the gaps among
the $m-1$ possibilities between $m+1$ and $2m-1$. Hence, there are
$\binom{m-1}{4}=\binom{\gamma -4}{4}$ options. If $\#A_{1}=3$, then
$\#A_{2}=1$ and, since $\max (A_{2})\leqslant 5$, there are only five options:
If $A_{2}=\{1\}$, then $A_{1}$ can be any subset of three elements between
$1$ and $m-1$ containing $1$. There are
$\binom{m-2}{2}=\binom{\gamma -5}{2}$ such options; If $A_{2}=\{2\}$, then
$A_{1}$ can be any subset of three elements containing $\{1,2\}$. There
are $m-3=\gamma -6$ options; If $A_{2}=\{3\}$, then $A_{1}$ can be any
subset of three elements either containing $\{1,3\}$ or containing
$\{2,3\}$. There are $2(m-3)-1=2\gamma -13$ options; If
$A_{2}=\{4\}$, then $A_{1}$ is either $\{1,2,4\}$ or $\{2,3,4\}$. Those
are $2$ options; If $A_{2}=\{5\}$, then $A_{1}$ is either
$\{1,2,5\},\{1,3,5\},\{2,4,5\},\{3,4,5\}$. Those are $4$ options. If
$\#A_{1}=2$, then $\#A_{2}=2$ and there are only two options: either
$A_{1}=A_{2}=\{1,2\}$ or $A_{1}=A_{2}=\{1,3\}$.

We conclude that
$n_{\gamma ,\gamma -3}=\binom{\gamma -4}{4}+\binom{\gamma -5}{2}+
\gamma -3\gamma -11$ and the result of the lemma follows.
\end{proof}

\section{The algorithm \textsc{ExploreUnleavedTree}}
\label{sec6}

All these results give rise to the trimming version of the algorithm
\Call{ExploreTree}{} when used to count semigroups of genus
$\gamma $, for $\gamma \geqslant 8$, without needing to visit them. It
is the algorithm \textsc{ExploreUnleavedTree} shown in \reftext{Algorithm~\ref{a:exploreunleavedtree}}.

\begin{algorithm}\caption{{\sc ExploreUnleavedTree}}\label{a:exploreunleavedtree}
{\small    \begin{algorithmic}
  \Function{ExploreUnleavedTree}{$\gamma$}
  \State{count $\gets\binom{\gamma-4}{4}+\binom{\gamma-2}{3}+\binom{\gamma-5}{2}+\lfloor\frac{\gamma}{3}\rfloor+6\gamma-12$}\Comment{{\color{blue}Lemma~\ref{l:minmultiplicities} and Lemma~\ref{l:maxmultiplicities}}}
    \For {$m$ \From $4$ \To $\gamma-4$}
  \State{$min\gets\min\{m,\gamma-m\}$}
 \For{$u$ \From 2 \To $min-1$}
 \State {count $\gets$ count$+$\Call{PseudoDescendAndTrim}{$m,{\mathbb N}_0,m+u,m,u,m-2,\gamma$}} \Comment{{\color{blue}$\omega(P_{m,u})=m$, $\shrink{P_{m,u}}={\mathbb N_0}$}}
 \EndFor
  \If{$min<\gamma-m$}
  \State {count $\gets$ count$+$\Call{PseudoDescendAndTrim}{$m, {\mathbb N}_0,2m,m,m,m-1,\gamma$}}  
  \Comment{{\color{blue}$\omega(P_{m,m})=m$, $\shrink{P_{m,m}}={\mathbb N_0}$}}
  \Else
  \State{count $\gets$ count$+\gamma-m+\binom{m-1}{2}$}
  \If{$2\gamma>3m$}
  \State{count $\gets$ count$-1$}
  \EndIf
  \EndIf
  \State{$r\gets m-3$}
  \For{$\sigma$ \From $m+2$ \To $2m-2$ {\bfseries while} genus$(\Lambda_{m,\dots,\sigma-1})\geq \gamma$} 
  \If{genus$(\Lambda_{m,\dots,\sigma-1})= \gamma$}
  \State{count $\gets$ count$+1$}
  \Else
  \State{count$\gets$count$+$\Call{DescendAndTrim}{$1,\Lambda_{\{m,\dots,\sigma-1\}},m,\sigma+1,m,1,r,\gamma$}}
  \Comment{{\color{blue}$\omega(Q_{m,\sigma})=1$, $\shrink{(Q_{m,\sigma})}=\Lambda_{\{m,\dots,\sigma-1\}}$}}
  \State{$r\gets r-1$}
  \EndIf
  \EndFor
  \EndFor
  \State{\Return count}
  \EndFunction
  \end{algorithmic}}
\end{algorithm}

We now want to have a look at the number of numerical semigroups that are
encoded in the algorithm \Call{ExploreUnleavedTree}{}. Note that there
are semigroups which are encoded and then trimmed after checking that they
have no descendants of genus $\gamma $. On the other hand, there are semigroups
which belong to the unleaved tree but need not to be encoded, as for instance
the semigroups of smallest or largest possible multiplicities, the pseudo-ordinary
semigroups, or the semigroups whose descendants in ${\mathscr{T}}$ have
maximum depth equal to $\gamma $, because we already know how many descendants
they have of genus~$\gamma $.

In \reftext{Fig.~\ref{f:tresarbresencoded}} we represented the complete tree of
genus $9$ (left), the subset of those semigroups that are encoded by the
algorithm (center), and the subset of those semigroups that are not encoded
(right).

\mut{
\begin{figure}
\includegraphics{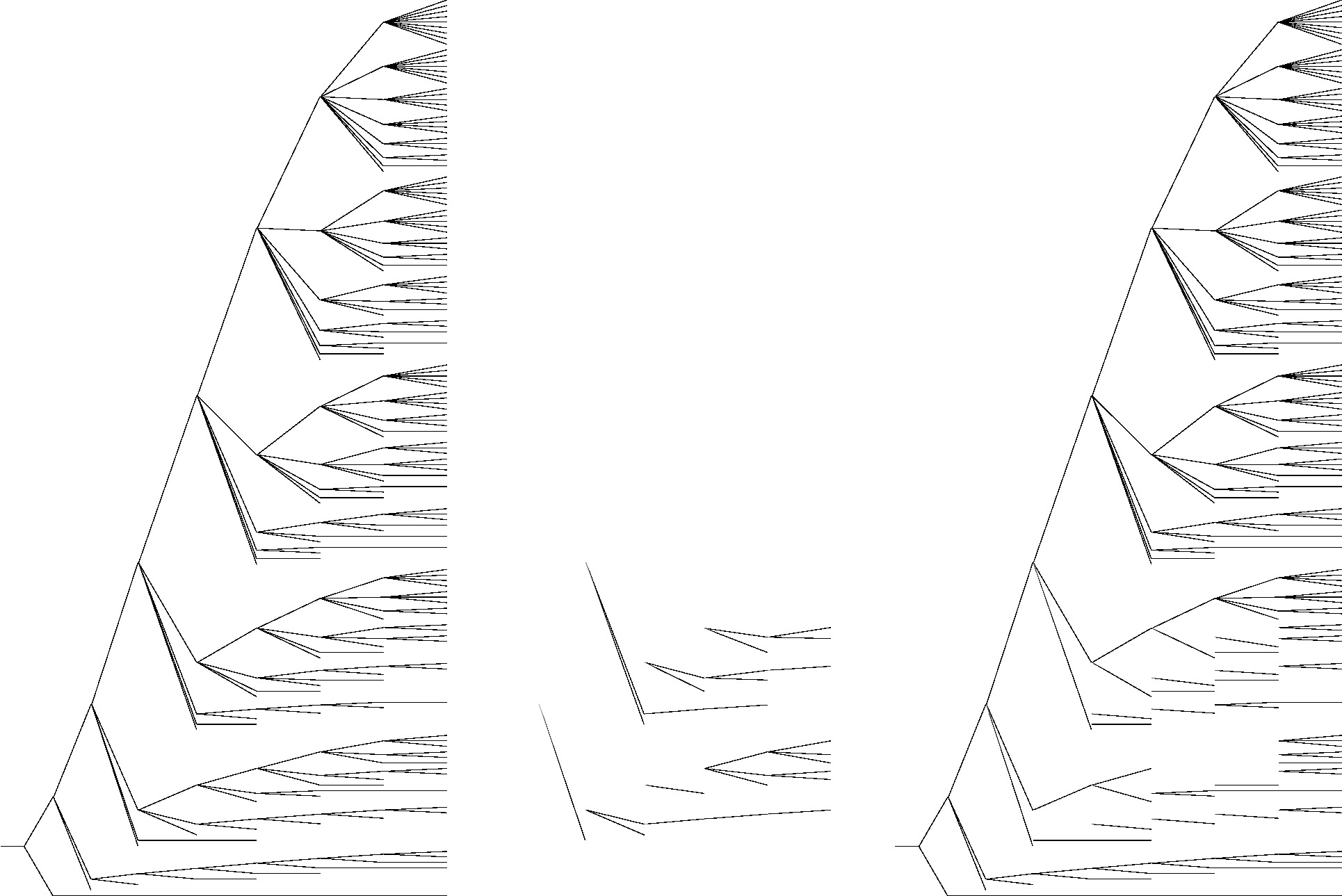}
\caption{All semigroups of genus up to $9$; subset of semigroups encoded
by the algorithm; subset of non encoded semigroups.}%
\label{f:tresarbresencoded}
\end{figure}
}

\begin{figure}
  \caption{All semigroups of genus up to $9$; subset of semigroups encoded by the algorithm; subset of non encoded semigroups}\label{f:tresarbresencoded}
    \resizebox{\textwidth}{!}{\includegraphics{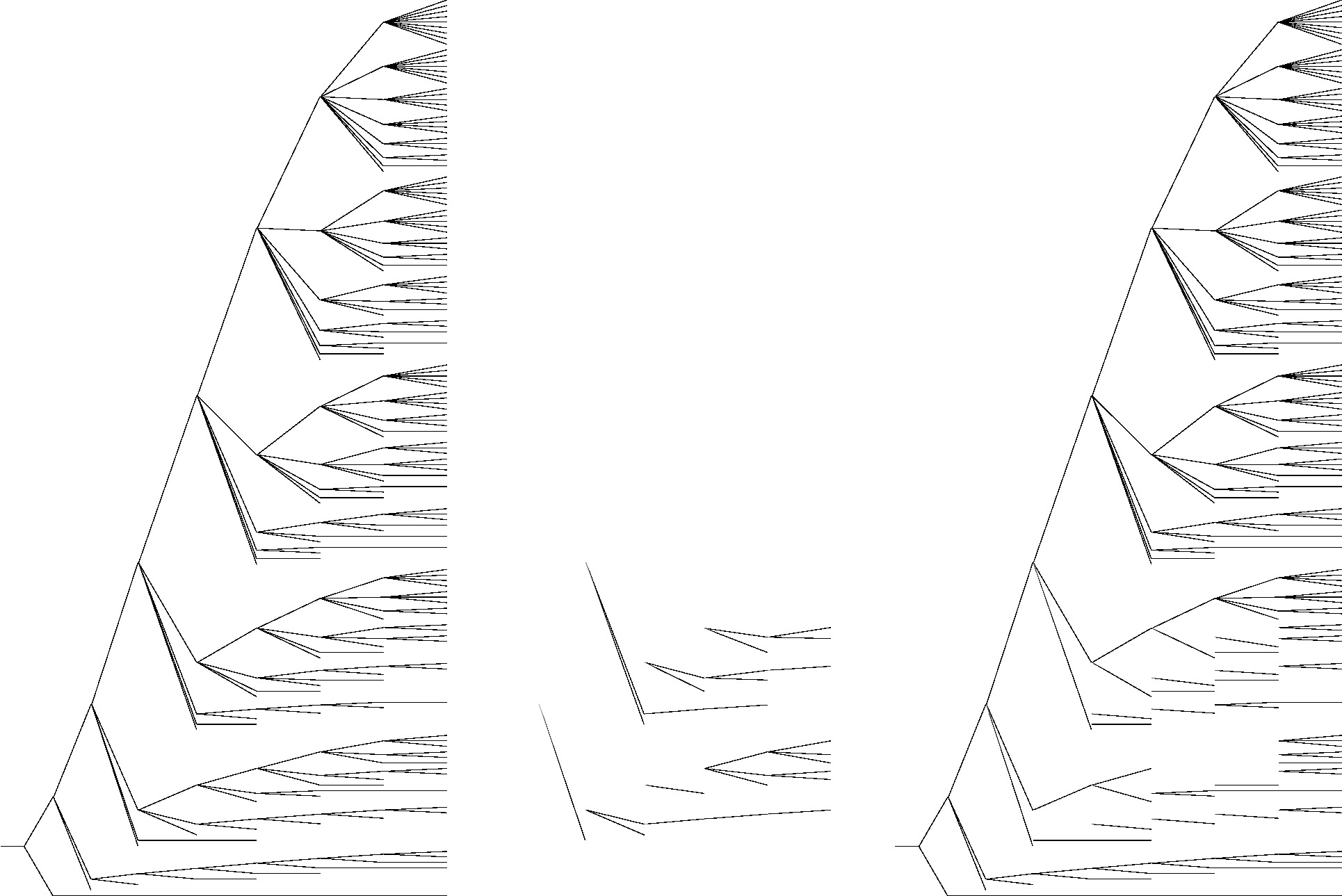}}
\end{figure}

\begin{table}
\caption{For each genus $\gamma $, number $n_{\gamma}$ of numerical semigroups
of genus $\gamma $, number of nodes in the complete tree of genus
$\gamma $, number of nodes in the unleaved tree of genus $\gamma $, number
of nodes encoded by the algorithm for genus $\gamma $, and their respective
percentages with respect to the number of nodes in the complete tree.}%
\label{t:t}
\resizebox{\textwidth}{!}{%
\begin{tabular}{llllllll}%
\hline $\gamma $& 10 & 15 & 20 & 25 & 30 & 35 & 40 \\\hline
  $n_\gamma $& 204 & 2857 & 37396 & 467224 & 5646773 & 66687201 & 774614284
  \\[\tvspace{6pt}]
 complete tree & 478 & 6964 & 93142 & 1179597 & 14396338 & 171202690 & 1998799015
 \\[\tvspace{6pt}]
 unleaved tree & 364  & 4833  & 61469   & 759972   & 9146174  & 107815637  & 1251716100  \\
  (\% wrt. complete)&  ($76\%$) &  ($69\%$) & ($66\%$) & ($64\%$) & ($64\%$) & ($63\%$) & ($63\%$)
  \\[\tvspace{6pt}]
  encoded nodes & 61  & 1325  & 16774  & 196433  & 2282567  & 26454236 & 304794995  \\
 (\% wrt. complete) &  ($13\%$) &  ($19\%$) &   ($18\%$) &  ($17\%$) & ($16\%$) & ($15\%$) & ($15\%$) \\\hline
\end{tabular}}%
\end{table}

The number of encoded nodes is expected to be proportional to the computation
time. In \reftext{Table~\ref{t:t}} we give the number of nodes in the complete tree,
number of nodes in the unleaved tree and the number of nodes encoded by
the algorithm \Call{ExploreUnleavedTree}{} for $\gamma $ from $10$ to
$40$. We can see that the number of nodes in the unleaved tree is around
$63\%$ of the number of nodes of the complete tree. Also, the number of
semigroups encoded by our algorithm is just around $15\%$ of the number
of nodes of the complete tree, and it is also less than half $n_{g}$. This
is the main reason for our algorithm to go much faster than any algorithm
exploring all nodes of ${\mathscr T}_{g}$.

\section*{Acknowledgement}

This research is supported by project PID2024-156636NB-C21 (MATSE), funded by MCIN/AEI/10.13039/501100011033/ FEDER, UE, and the project RED2024-153572-T, funded by MICIU/AEI/10.13039/501100011033.

The author is grateful to Christopher O'Neill for very interesting discussions and for suggesting the formula for $n_{\gamma,\gamma-3}$ used in the proof of Lemma~\ref{l:maxmultiplicities}. She is also grateful to Jose Luis Muñoz, Armando Martín López and Juan Paz Sánchez for facilitating the computational resources that made it possible to compute $n_{76}$ and $n_{77}$.


\end{document}